\documentclass[8pt,draftcls,twocolumn]{autart}

\usepackage{graphicx}          % Include this line if your
\usepackage{amsfonts}
\usepackage{epsfig}
\usepackage{subfigure}
\usepackage{amsmath,color}
\usepackage[longnamesfirst,round]{natbib}
\usepackage[all]{xy}
\usepackage{float}
\usepackage{mathrsfs}
\usepackage{amssymb}
\usepackage{CJK}

\newtheorem{definition}{Definition}[section]

\newtheorem{theorem}{Theorem}[section]

\newtheorem{remark}{Remark}[section]
\newtheorem{assumption}{Assumption}[section]

\def\proof{\noindent{\it Proof: }}
\def\QED{\mbox{\rule[0pt]{1.5ex}{1.5ex}}}
\def\endproof{\hspace*{\fill}~\QED\par\endtrivlist\unskip}

\def\T{^{\rm\tiny T}}
\begin{document}
\begin{frontmatter}
\title{Behaviors of Networks with Antagonistic Interactions and Switching Topologies \thanksref{footnoteinfo}}

\thanks[footnoteinfo]{This paper was not presented at any IFAC
meeting. Corresponding author: Z. Meng.}

\author[Tsinghua]{Ziyang Meng}\ead{ziyangmeng@mail.tsinghua.edu.cn},
\author[ANU]{Guodong Shi}\ead{guodong.shi@anu.edu.au},
\author[KTH]{Karl H. Johansson}\ead{kallej@kth.se},
\author[RUG]{Ming Cao}\ead{m.cao@rug.nl},
\author[CAS]{Yiguang Hong}\ead{yghong@iss.ac.cn}

\address[Tsinghua]{State Key Laboratory of Precision Measurement Technology and Instruments, Department of Precision Instrument, Tsinghua University, Beijing 100084, China.}
\address[ANU]{College of Engineering and Computer Science, The Australian National University, Canberra, Australia.}
\address[KTH]{ACCESS Linnaeus Centre and School of Electrical Engineering,
Royal Institute of Technology, Stockholm, Sweden.}
\address[RUG]{Engineering and Technology Institute, University of Groningen, the Netherlands.}
\address[CAS]{Key Laboratory of Systems and Control, Institute of Systems Science, Chinese Academy of Sciences, Beijing, China.}

%%%%%%%%%%%%%%%%%%%%%%%%%%%%%%%%%%%%%%%%%%%%%%%%%%%%%%%%%%%%%%%%%%%%%%%
%\begin{keyword}
%Network Behavior; Antagonistic Interactions; Switching Topologies
%\end{keyword}

\begin{abstract}
In this paper, we study the discrete-time consensus problem over networks with antagonistic and cooperative interactions.
A cooperative interaction between two nodes takes place when one node receives the true state of the other while
an antagonistic interaction happens when the former receives the opposite of the true state of the latter. We adopt a quite general model where the node communications can be either unidirectional or bidirectional, the network topology graph may vary over time, and the cooperative or antagonistic relations can be time-varying. It is proven that, the limits of all the node states exist, and the absolute values of the  node states reach consensus if the switching interaction graph is uniformly jointly strongly connected for unidirectional topologies, or infinitely jointly connected for bidirectional topologies. These results are independent of the switching of the interaction relations. We construct a counterexample to indicate a rather surprising fact that quasi-strong connectivity of the interaction graph, i.e., the graph contains a directed spanning tree, is not sufficient to guarantee the consensus in absolute values even under fixed topologies. Based on these results, we also propose sufficient conditions for bipartite consensus to be achieved over the network with joint connectivity. Finally, simulation results using a discrete-time Kuramoto model are given to illustrate the convergence results showing that the proposed framework is applicable to a class of networks with general nonlinear dynamics.
\end{abstract}

\end{frontmatter}

%%%%%%%%%%%%%%%%%%%%%%%%%%%%%%%%%%%%%%%%%%%%%%%%%%%%%%%%%%%%%%%%%%%%%%%%%
\section{Introduction}
%Consensus seeking over multi-agent networks has been extensively studied during the past decade due to its potential applications in various areas including spacecraft formation flying (\cite{WangHadaegh96}), control of multiple unmanned aerial vehicles (\cite{LinBrouckeFrancis04}), distributed estimation of sensor networks (\cite{Schenato_IJSC08}), and collective behaviors of biological swarming (\cite{TannerJadbabaiePappas07}).
%In particular,
Distributed consensus algorithms were first introduced in the study of distributed optimization methods in (\cite{Tsitsiklis_TAC86}). Phase synchronization was observed in (\cite{VicsekEtAl95}) and its mathematical proof was given in (\cite{JadbabaieLinMorse03}). The robustness of the consensus algorithm to link/node failures and time-delays was studied in (\cite{SaberFaxMurray07_IEEE}).
A central problem in consensus study is to investigate the influence of the interaction graph on the convergence or convergence speed of the multi-agent system dynamics. The interaction graph, which describes the information flow among the nodes, is often time-varying due to the complexity of the
interaction patterns in practice. Both continuous-time and discrete-time models were studied for consensus algorithms with switching interaction graphs and joint connectivity conditions were established for linear models (\cite{Blondel_CDC2005,RenBeard05_TAC,CaoMing_SIAM08J2,HendrickxTAC13}). Nonlinear multi-agent dynamics have also drawn much attention (\cite{ShiGuodong_Automatica2009, MengZiyang_Automatica2013}) since in many practical problems the node dynamics are naturally nonlinear, e.g., the Kuramoto model (\cite{kuramoto1}).
%It was shown in (\cite{Moreau_TAC05}) that convexity, instead of linearity, is essentially important in reaching a consensus.

Although great progress has been made, most of the existing results are based on the assumption that the agents in the network are cooperative. Recently, motivated by opinion dynamics over social networks (\cite{Harary1956,Krause2002,Kleinberg_Book}),
consensus algorithms over cooperative-antagonistic networks drew much attention (\cite{AltafiniPlosOne,Altafini_TAC13,Shi_JSAC2013}). Altafini (\cite{Altafini_TAC13}) assumed that a node receives the opposite of the true state of its neighboring node if they are antagonistic. Therefore, the modeling of such an antagonistic input for agent $i$ is of the form $-(x_i+x_j)$ (in contrast to the form $-(x_i-x_j)$ for cooperative input), where $j$ denotes the antagonistic neighbor of agent $i$.
On the other hand, the authors of \cite{Shi_JSAC2013} assumed that a node receives the opposite of the relative state from its neighboring node if they are antagonistic. Then, the antagonistic input for agent $i$ is modeled by the form $(x_i-x_j)$ in this case. The extension to the case of homogenous single-input high-order dynamical systems was discussed in (\cite{Valcher_SCL14}). The graph was assumed to be fixed and a spectral analysis approach was used. Instead, we will focus on switching topologies with joint connectivity and take advantage of a detailed state-space analysis approach in this paper.
A lifting approach was proposed in \cite{Hendrickx_CDC14} to study opinion dynamics with antagonisms over switching interaction graphs. Some general conditions were established by applying the rich results from the consensus literature. The dissensus problem was studied in (\cite{Bauso_TAC12}), where the focus was to understand when consensus is or is not achieved if death and duplication phenomena occur.
Note that in dissensus the control terms for death and
duplication phenomena are added to the classical consensus algorithm. However, consensus or bipartite
consensus in our study does not denote a control goal, but a final behavior of multi-agent systems. In particular, consensus denotes the final states of all the agents converging to the same value while bipartite
consensus denotes the final states of all the agents converging to two
opposite values.

Note that most of the existing works on antagonistic interactions are based on the assumption that the interaction graph is fixed. In many practical cases, however, the interactions between agents may vary over time or be dependent on the states. In this paper, we focus on the behavior of multiple agents with antagonistic interactions, discrete-time dynamics, and switching interaction graphs. Both unidirectional and bidirectional topologies are considered. We show that the limits of all node states exist and reach a consensus in absolute values if the switching interaction graph is uniformly jointly strongly connected for unidirectional topologies, or infinitely jointly connected for bidirectional topologies. Here, reaching consensus in absolute values is not a design objective, but rather an emergent behavior that
we can observe for cooperative-antagonistic networks.
By noting that an antagonistic interaction also represents an arc in the graph, we know that the connectivity of the network may not be guaranteed without antagonistic interactions. Therefore, we actually show that the antagonistic interaction has a similar role as the cooperative interaction in contributing to the consensus of the absolute values of the node states. In addition,
a counterexample is constructed that indicates that quasi-strong connectivity of the interaction graph, i.e., the graph has a directed spanning tree, is not sufficient to guarantee consensus of the node states in absolute
values even under a fixed topology. Based on these results, we propose sufficient conditions for bipartite consensus to
be achieved over a network with joint connectivity. It turns out that the structural balance condition is essentially important and this part of the result can be viewed as an extension of the work \cite{Altafini_TAC13} to the case of
general time-varying graphs with joint connectivity. A detailed
asymptotic analysis is performed with a contradiction argument to show the main results.
%
%The remainder of the paper is organized as follows. In
%Section \ref{sec:ps}, we give the problem formulation and present the  main results. The  proofs of the results are given in Section \ref{sec:proofs}.  A simulation example and some brief concluding remarks are given in Sections \ref{sec:simulation}
%and \ref{sec:conclusion}, respectively.

\section{Problem Formulation and Main Results} \label{sec:ps}
Consider a multi-agent network with agent set $\mathcal{V}=\{1,\dots,n\}$. In the rest of the paper we use  {\it agent} and {\it node} interchangeably. The state-space for the agents is $\mathbb{R}$, and we let $x_i\in\mathbb{R}$ denote the state of node $i$. Set $x=(x_1, x_2,\dots, x_n)\T$.

\subsection{Interaction Graph}\label{sec:graph}
The interaction graph of the network is defined as a sequence of unidirectional graphs, $\mathcal{G}_{k}=(\mathcal{V}, \mathcal{E}_{k})$, $k=0,1,\dots$, with node set $\mathcal{V}$ and $\mathcal{E}_k \subseteq \mathcal{V}\times \mathcal{V}$ is the set of arcs at time $k$. An arc from node $i$ to $j$ is denoted by $(i,j)$. A path from node $i$ to $j$ is a sequence of consecutive arcs $ \{(i, k_1),\dots, (k_l, j)\}$. We assume that $\mathcal{G}_{k}$ is a signed graph, where ``+'' or ``$-$'' is associated with each arc $(i,j)\in \mathcal{E}_k$. Here, ``+'' represents cooperative relation and ``$-$'' represents antagonistic relation. The set of neighbors of node $i$ in $\mathcal{G}_{k}$ is denoted by $\mathcal{N}_i(k):=\{j:(j,i)\in \mathcal{E}_k\}\cup \{i\}$, and $\mathcal{N}_i^+(k)$ and $\mathcal{N}^-_i(k)$ are used to denote the cooperative neighbor sets and antagonistic neighbor sets, respectively. Clearly, $\mathcal{N}_i(k)=\mathcal{N}_i^-(k)\cup\mathcal{N}_i^+(k)\cup \{i\}$.
The joint graph of $\mathcal{G}$ during time interval $[k_1,k_2)$ is defined by $\mathcal{G}([k_1,k_2))=\bigcup_{k\in[k_1,k_2)}\mathcal{G}_k=(\mathcal{V},
\bigcup_{k\in[k_1,k_2)}\mathcal{E}_{k})$. The sequence of graphs $\{\mathcal{G}_{k}\}_0^\infty$ is said to be sign consistent if the sign of any arc $(i,j)$ does not change over time. Under the assumption that $\{\mathcal{G}_{k}\}_0^\infty$ is sign consistent, we can define a signed total graph $\mathcal{G}^*:=(\mathcal{V}, \mathcal{E}^*)$, where $\mathcal{E}^*=\cup_{k=0}^{\infty}\mathcal{E}_k$.

We write $i \rightarrow j$ if there is a path from $i$ to $j$. A root is a node $i$ such that  $i \rightarrow j$ for every other node $j\in\mathcal{V}\backslash\{i\}$. A unidirectional graph is quasi-strongly connected if it has a directed spanning tree, i.e., there exists at least one root.
%We say node $j$ is reachable from node $i$ in a digraph if there exists a path from $i$ to $j$. A unidirectional graph is quasi-strongly connected if it has a directed spanning tree, i.e., there exists at least one node that is reachable to all the  other nodes.
A unidirectional graph is called strongly connected if there is a path connecting any two distinct nodes.
A unidirectional graph $\mathcal{G}$ is called bidirectional if for any two nodes $i$ and $j$, $(j,i)\in\mathcal{E}$ if and only if $(i,j)\in\mathcal{E}$. A bidirectional graph is connected if it is connected as a bidirectional graph ignoring the arc directions. We introduce the following definition of the joint connectivity of a sequence of graphs.

\begin{definition} (i). $\{\mathcal{G}_{k}\}_0^\infty$ is uniformly jointly strongly connected if there exists a constant $T\geq 1$ such that $\mathcal{G}([k,k+T))$ is strongly connected for any $k\geq 0$.

(ii). $\{\mathcal{G}_{k}\}_0^\infty$ is uniformly jointly quasi-strongly connected if there exists a constant $T\geq 1$ such that $\mathcal{G}([k,k+T))$ has a directed spanning tree for any $k\geq 0$.

(iii). Suppose $\mathcal{G}_{k}$ is bidirectional for all $k\geq0$. Then $\{\mathcal{G}_{k}\}_0^\infty$ is infinitely jointly connected if $\mathcal{G}([k,+\infty))$ is connected for any $k\geq 0$.
\end{definition}

\subsection{Node Dynamics}\label{sec:control}

The update rule for each node is described by:
\begin{align}
x_i(k+1)=&\sum_{j\in \mathcal{N}_i(k)}a_{ij}(x,k)x_j(k),\notag \\& k=0,1,\dots,
\quad i=1,2,\dots,n,\label{eq:origin}
\end{align}
where $x_i(k)\in\mathbb{R}$ represents the state of agent $i$ at time $k$, $x=(x_1,x_2,\dots,x_n)\T$, and $a_{ij}(x,k)$ represents a nonlinear time-varying function.
%the weight of arc $(j,i)$ associated with the interaction graph $\mathcal{G}_{k}$.
Equation \eqref{eq:origin} can be written in the compact form:
\begin{equation}
x(k+1)=A(x,k)x(k),\quad k=0,1,\dots,\label{eq:matrix}
\end{equation}
where $A(x,k)=[a_{ij}(x,k)]\in \mathbb{R}^{n\times n}$.
For $a_{ij}(x,k)$, we impose the following assumption.
\vspace{2mm}

\begin{assumption}\label{assm1}
There exists a positive constant $0<\lambda< 1$ such that:

(i) $a_{ii}(x,k)\geq \lambda$, for all $i,x,k$, and $\sum_{j\in \mathcal{N}_i(k)}|a_{ij}(x,k)|=1$ for all $i,x,k$;

(ii) $a_{ij}(x,k)\geq \lambda$ for all $i,j,x,k$, if $j\in \mathcal{N}_i^+(k)$; and $a_{ij}(x,k)\leq -\lambda$ for all $i,j,x,k$, if $j\in \mathcal{N}_i^-(k)$.
\end{assumption}

\begin{remark}
The antagonistic interactions commonly exist in social networks and signed graphs are used to describe these interactions. In signed graphs, a positive/negative weight is associated with a cooperative/antagonistic relationship between
the two agents.
Assumption \ref{assm1} models such a signed graph in a discrete-time dynamics and switching topology setting. Similar modeling can be found in \cite{Altafini_TAC13} with continuous-time dynamics and fixed topology.
%Firstly, Assumption \ref{assm1} is motivated by the definition of signed graphs with continuous-time dynamics and fixed topology (see Section II of \cite{Altafini_TAC13}).
Clearly, the second part of Assumption \ref{assm1} describes cooperative-antagonistic interactions in the sense that $a_{ij}(x,k)>0$ represents that $i$ is cooperative to $j$, and $a_{ij}(x,k)<0$ represents that $i$ is antagonistic to $j$.
%If $a_{ij}$ is assumed to be constant, Assumption \ref{assm1} is reduced to the discrete-time version of the model studied in %\cite{Altafini_TAC13} since the existence of $\lambda$ is trivially satisfied in such a case.
In addition, because of complexity of the
interaction patterns, $a_{ij}$ may depend on time or relative measurements for agent dynamics \eqref{eq:matrix}, instead of being constant.
Lots of practical multi-agent system models can be written in this form (e.g., Kuramoto equation, consensus algorithm, and swarming model given in \cite{Moreau_TAC05}). Last but not the least, we want to emphasize that
the existence of a positive constant $\lambda$ in Assumption \ref{assm1} is a technical assumption and plays an indispensable role in driving the states of the system asymptotically to
converge (see the proofs of the main theorems). This is in fact a very general assumption and has been used extensively in the existing literature (see e.g., \cite{Blondel_CDC2005,Shi_JSAC2013}). Here, $\lambda$ can be an arbitrarily small constant and thus the assumption on the existence of $\lambda$ will not restrict the usefulness of Assumption \ref{assm1} on practical applications.
%thus the interaction can be very weak (small $a_{ij}$) possible.
%In fact, when $a_{ij}$ is chosen from a finite set (e.g., the case of linear time-varying dynamics), the existence of such a $\lambda$ is trivially satisfied.

%Thirdly, $a_{ij}(x,k)$ does not depend on global information in \eqref{eq:matrix}. For example, $a_{ij}(x,k)$ can only depend on the state of $x_i$, $x_j$, $j\in\mathcal{N}_i$ and time $k$, which maintain the spirit of distributed control.
\end{remark}

\subsection{Main Results}
Uniform joint strong connectivity is sufficient for convergence of agent states for unidirectional graphs, as stated in the following theorem.
\begin{theorem}\label{thm1}
Suppose that Assumption \ref{assm1} holds and that $\mathcal{G}_{k}$ is unidirectional for all $k\geq 0$. For system \eqref{eq:origin}, $\lim_{k\rightarrow\infty}x_i(k)$ exists and $\lim_{k\rightarrow\infty}|x_i(k)|=M^*$, for all $i\in\mathcal{V}$ and
every initial state $x(0)\in \mathbb{R}^{n}$, if $\{\mathcal{G}_{k}\}_0^\infty$ is uniformly jointly strongly connected, where $M^*$ is a nonnegative constant.
\end{theorem}

%\begin{remark}
%It has been shown that consensus may not be achieved
%for cooperative-antagonistic multi-agent system; see Example 1 of \cite{Altafini_TAC13}. Instead, Altafini showed that bipartite consensus can be achieved if the sign-symmetric signed graph is strongly connected and structurally balanced.
%Compared to the results given in (\cite{Altafini_TAC13}), Theorem \ref{thm1} requires no conditions on
%structurally balance of the sign graph $\mathcal{G}_k$. By noting that an antagonistic interaction also represents an arc in the graph, we know that the connectivity of the network may not be guaranteed without antagonistic interactions. According to the proofs of Theorem \ref{thm1}, the connectivity of the network is indispensable to derive the conclusions and, therefore, we have actually shown that both cooperative and antagonistic interactions contribute to the consensus of the absolute values of the node states.
%\end{remark}

For cooperative networks, it is well-known that asymptotic consensus can be achieved if the interaction graph is uniformly quasi-strongly connected, e.g., (\cite{RenBeard05_TAC,CaoMing_SIAM08J1}). Note that quasi-strong connectedness is weaker than strong connectedness. Then, a natural question is whether asymptotic consensus in absolute values can be achieved when the interaction graph is uniformly quasi-strongly connected. We construct the following counterexample showing that quasi-strong connectivity is not sufficient for the consensus in absolute values even in case of a fixed graph.

\vspace{2mm}

\noindent {\bf Counterexample.} Let $\mathcal{V}=\{1,2,3\}$ and initial state is $x(0)=(1,0,-1)\T$. The interaction graph $\mathcal{G}_k=\mathcal{G}$ is fixed and shown in Fig. 1 and
\begin{displaymath}
A=[a_{ij}]=\left[\begin{matrix} 1&0&0\\1/3 &1/3& 1/3\\-1/2 &0 &1/2
\end{matrix}\right].
\end{displaymath}
It is straightforward to check that $\mathcal{G}$ is quasi-strongly connected. However, the states of the agents remain $x_1(k)=1$, $x_2(k)=0$, and $x_3(k)=-1$, for all $k=1,2,\dots$ under system dynamics \eqref{eq:origin}. Thus, the absolute values of the agent states do not reach a consensus.

\begin{figure}[t]\centering
\begin{tabular} {c}
 \xymatrix{
1 \ar@{->}[rd]_-{+}\ar@{->}[r]_-{-} &  3\ar@{->}[d]_-{+}\\
&  2
} \end{tabular}\caption{Signed graph $\mathcal{G}$ used in the counterexample} \label{fig1}
\end{figure}
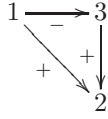

For bidirectional graphs, we present the following result indicating that the consensus in absolute values can be achieved under weaker connectivity conditions than those in Theorem \ref{thm1} for unidirectional graphs.

\begin{theorem}\label{thm2}
Suppose that Assumption \ref{assm1} holds and $\mathcal{G}_{k}$ is bidirectional for all $k\geq 0$. For system \eqref{eq:origin}, $\lim_{k\rightarrow\infty}x_i(k)$ exists and $\lim_{k\rightarrow\infty}|x_i(k)|=M^*$, for all $i\in\mathcal{V}$ and every initial state $x(0)\in \mathbb{R}^{n}$ if $\{\mathcal{G}_{k}\}_0^\infty$ is infinitely jointly connected, where $M^*$ is a nonnegative constant.
\end{theorem}

We note that Theorems \ref{thm1} and \ref{thm2} are concerned with consensus of the absolute values of agent states. It is not clear what is the final sign of the states, i.e., the sign of $\lim_{k\rightarrow\infty}x_i(k)$ for different $i\in\mathcal{V}$. We next characterize how bipartite consensus, i.e., splitting into two opposite states, can emerge. We first introduce the notion of structural balance, cf., Definition 2 of \cite{Altafini_TAC13}.
\begin{definition}
Suppose that the sequence of graphs $\{\mathcal{G}_{k}\}_0^\infty$ is sign consistent and $\mathcal{G}^*=(\mathcal{V},\mathcal{E}^*)$ is a signed total graph defined in Section \ref{sec:graph}. $\mathcal{G}^*$ is structurally balanced if we can divide $\mathcal{V}$ into two disjoint nonempty subsets $\mathcal{V}_1$ and $\mathcal{V}_2$ (i.e., $\mathcal{V}_1\bigcup\mathcal{V}_2=\mathcal{V}$ and $\mathcal{V}_1\bigcap\mathcal{V}_2=\emptyset$), where negative arcs only exist between these two subsets, i.e., the arc $(i,j)$ is associated with sign ``+'', $\forall i,j\in\mathcal{V}_q$ ($q\in\{1,2\}$) and the arc $(i,j)$ is associated with sign ``$-$'', $\forall i\in\mathcal{V}_q,j\in\mathcal{V}_r,q\neq r$ ($q,r\in\{1,2\}$).
\end{definition}
\begin{remark}
Since we assume that the sequence of graphs $\{\mathcal{G}_{k}\}_0^\infty$ is sign consistent, even if the interconnection graph may be changing with time,
the sets ${\mathcal V}_1$ and ${\mathcal V}_2$ do not, and hence the two ``antagonistic groups''
are always the same given that the signed total graph $\mathcal{G}^*$ is structurally balanced.
\end{remark}
We next establish our result for bipartite consensus.
\begin{theorem}\label{thm3}
 Suppose that $\mathcal{G}_{k}$ is unidirectional for all $k\geq 0$ and $\{\mathcal{G}_{k}\}_0^\infty$ is uniformly jointly strongly connected, or $\mathcal{G}_{k}$ is bidirectional for all $k\geq0$ and $\{\mathcal{G}_{k}\}_0^\infty$ is infinitely jointly connected. Also suppose that Assumption \ref{assm1} holds,  $\mathcal{G}^*$ is structurally balanced, and every negative arc in $\mathcal{G}^*$ appears infinitely often in $\{\mathcal{G}_{k}\}_0^\infty$. For system \eqref{eq:origin}, $\lim_{k\rightarrow\infty}x_i(k)=M^*$, $i\in\mathcal{V}_1$, and $\lim_{k\rightarrow\infty}x_i(k)=-M^*$, $i\in\mathcal{V}_2$, for every initial state $x(0)\in \mathbb{R}^{n}$ , where $M^*$ is a constant.
\end{theorem}
\begin{remark}
The interpretation of Theorem \ref{thm3} is that under proper graph structure, the agent states in the same cooperative subgroup  converge to the same limit and the limits of agents in different cooperative subgroups are exactly opposite.
It is not hard to see from Theorem \ref{thm3} that the structural balance condition plays a key role in obtaining such a bipartite consensus behavior. Theorem~\ref{thm3} can be viewed as an extension of the result of \cite{Altafini_TAC13} from fixed to switching graphs. We establish conditions for both unidirectional and bidirectional graphs with joint connectivity. A detailed analysis of the asymptotic behavior will be done and a contradiction argument applied to show the main results, in contrast to the spectral analysis approach given in \cite{Altafini_TAC13}.
\end{remark}
%\begin{theorem}\label{thm4}
%Suppose that Assumption \ref{assm1} holds, $\mathcal{G}_{k}$ is bidirectional for all $k\geq0$ and $\{\mathcal{G}_{k}\}_0^\infty$ is infinitely jointly connected.
%Also suppose that $\mathcal{G}^*$ is structurally balanced and every negative arc in $\mathcal{G}^*$ appears infinitely often in $\{\mathcal{G}_{k}\}_0^\infty$. For system \eqref{eq:origin}, $\lim_{k\rightarrow\infty}x_i(k)=M^*$, $i\in\mathcal{V}_1$ and $\lim_{k\rightarrow\infty}x_i(k)=-M^*$, $i\in\mathcal{V}_2$ for every initial state $x(0)\in \mathbb{R}^{n}$ , where $M^*$ is a constant.
%\end{theorem}

\section{Proofs}\label{sec:proofs}
In this section, we present the proofs of the statements. First a key technical lemma is established, and then the proofs of Theorems \ref{thm1}, \ref{thm2}, and \ref{thm3} are presented. Since the proofs do not rely on whether or not $a_{ij}$ depends on $x$, without loss of generality, we use $a_{ij}(k)$ to denote $a_{ij}(x,k)$.

\begin{lem}\label{lem:invariant}
Suppose that Assumption \ref{assm1} holds. For system \eqref{eq:origin}, it holds that $\|x(k+1)\|_{\infty}\leq \|x(k)\|_{\infty}$, for all $k=0,1,\dots$.
\end{lem}
\proof
It follows from Assumption \ref{assm1} that
$
|x_{i}(k+1)|\leq ~\sum_{j\in \mathcal{N}_i(k)}|a_{ij}(k)||x_j(k)|
\leq  ~\left(\sum_{j\in \mathcal{N}_i(k)} |a_{ij}(k)|\right)
\\ \times\max_{i\in \mathcal{V}}|x_i(k)|
=~ \|x(k)\|_{\infty},
$
for all $i$, which leads to the conclusion directly.
\endproof

\subsection{Proof of Theorem \ref{thm1}: Consensus in Absolute Values }\label{sec:directed1}

Since a bounded monotone sequence always admits a limit,
Lemma \ref{lem:invariant} implies that for any initial value $x(0)$, there exists a constant $M^*$, such that
$\lim_{k\rightarrow \infty}\|x(k)\|_{\infty}=M^*$. We further define
$
\Phi_{i}=\limsup_{k\rightarrow \infty} |x_{i}(k)|, ~\Psi_{i}=\liminf_{k\rightarrow \infty} |x_{i}(k)|,
$
for all $i\in\mathcal{V}$.
Clearly, it must hold that $0\leq \Psi_{i}\leq \Phi_{i}\leq M^*$. Therefore, $\lim_{k\rightarrow \infty}|x_i(k)|=M^*$, for all $i\in\mathcal{V}$ if and only if $\Phi_{i}=\Psi_{i}=M^*$, $i\in \mathcal{V}$.
%We first use contradiction argument to show that $\lim_{k\rightarrow\infty}|x_i(k)|=M^*$, for all $i\in\mathcal{V}$.

In addition, based on the fact that $\lim_{k\rightarrow \infty}\|x(k)\|_{\infty}= M^*$, it follows that for any $\varepsilon>0$, there exists a $\widehat{k}(\varepsilon)>0$ such that $M^*-\varepsilon\leq\|x(k)\|_{\infty}\leq M^*+\varepsilon$, $\forall k\geq \widehat{k}(\varepsilon)$ and
$|x_{i}(k)|\leq M^*+\varepsilon,~~\forall i\in \mathcal{V}, ~~\forall k\geq \widehat{k}(\varepsilon).$

We next use a contradiction argument.
Now suppose that there exists a node $ {i_1}\in \mathcal{V}$ such that $0\leq \Psi_{i_1}<  M^*$. With the definitions of $\Psi_{i_1}$, for any $\varepsilon>0$, there exists a constant $\Psi_{i_1}<\alpha_{1}<M^*$
%and a subsequence $|x_i(k_0^*)|$, $|x_i(k_1^*)|$, $\dots$, $|x_i(k_l^*)|$, $\dots$, where $0\leq k_0^*<k_1^*<\dots<k_l^*<\dots$ is an increasing sequence of indices such that $|x_i(k_l^*)|\leq\alpha_{1}$, for all $l=0,1,\dots$. Therefore, it follows that there exists
and a time instance $k_1\geq \widehat{k}(\varepsilon)$ such that $|x_{i_1}(k_1)|\leq\alpha_{1}$. This shows that
\begin{align}
|x_{i_1}(k_1)|\leq M^*-(M^*-\alpha_{1}) = M^*-\xi_1,\label{eq:i1}
\end{align}
where $\xi_1=M^*-\alpha_{1}>0$.

First of all, it follows from Lemma \ref{lem:invariant} that $|x_{j}(k_1+s)|\leq \|x(k_1)\|_{\infty}$, for all $s=1,2,\dots$ and all $j\in\mathcal{V}$. We next fix $i_1$ and analyze the trajectory of $x_{i_1}$ after $k_1$.
Then it must be true that for all $s=1,2,\dots$,
$
|x_{i_1}(k_1+s)|\leq ~\sum_{j\in \mathcal{N}_{i_1}(k_1+s-1)}\left|a_{i_1j}(k_1+s-1)\right| \left|x_j(k_1+s-1)\right|
\\=~\left|a_{i_1i_1}(k_1+s-1)\right| \left|x_{i_1}(k_1+s-1)\right|
\\+\sum_{j\in \mathcal{N}_{i_1}(k_1+s-1)\setminus\{ i_1\} }
|a_{i_1j}(k_1+s-1)|
 |x_j(k_1+s-1)|
 \leq~|a_{i_1i_1}(k_1+s-1)||x_{i_1}(k_1+s-1)|
+(1-|a_{i_1i_1}(k_1+s-1)|)\|x(k_1)\|_{\infty}.
$
Also note that $\|x(k_1)\|_{\infty}\leq M^*+\varepsilon$.
It thus follows from \eqref{eq:i1} that
\begin{align}
&|x_{i_1}(k_1+1)|\nonumber
\\& \leq~ |a_{i_1i_1}(k_1)||x_{i_1}(k_1)|
+(1-|a_{i_1i_1}(k_1)|)\|x(k_1)\|_{\infty}\nonumber
\\& \leq ~a_{i_1i_1}(k_1)(M^*-\xi_1)
+(1-a_{i_1i_1}(k_1))(M^*+\varepsilon)\nonumber
\\& \leq ~M^*+\varepsilon-\lambda\xi_1,
\end{align}
where we have used the fact that $a_{i_1i_1}\geq \lambda$ from Assumption \ref{assm1}.

By a recursive analysis we can further deduce that
\begin{align}\label{r1}
|x_{i_1}(k_1+s)|\leq M^*+\varepsilon-\lambda^s\xi_1, \quad s=1,2,\dots.
\end{align}

Next, we consider the time interval $[k_1,k_1+T)$.
Since $\mathcal{G}([k_1,k_1+T))$ is strongly connected, there is a path from $i_1$ to any other node during the time interval $[k_1,k_1+T)$. This implies that
there exists a time instant $k_2\in [k_1,k_1+T)$ such that $i_1$ is a neighbor of another node $i_2$ at $k_2$. We next analyze the trajectory of $x_{i_2}$ after $k_2$. It follows that
$
|x_{i_2}(k_2+s)|\leq ~\sum_{j\in \mathcal{N}_{i_2}(k_2+s-1)}|a_{i_2j}(k_2+s-1)| |x_j(k_2+s-1)|
=~|a_{i_2i_1}(k_2+s-1)||x_{i_1}(k_2+s-1)|
+\sum_{j\in \mathcal{N}_{i_2}(k_2+s-1)\setminus\{i_1\}}\!\!|a_{i_2j}(k_2+s-1)| |x_j(k_2+s-1)|
\leq~|a_{i_2i_1}(k_2+s-1)| |x_{i_1}(k_2+s-1)|
+(1-|a_{i_2i_1}(k_2+s-1)|)\|x(k_2)\|_{\infty},
$
where we have used the fact that $ |x_i(k_2+s)|\leq \|x(k_2)\|_{\infty}$, for all $s=1,2,\dots$, and for all $i\in \mathcal{V}$.
Noting that $\|x(k_2)\|_{\infty}\leq M^*+\varepsilon$, it thus follows that
$
|x_{i_2}(k_2+s)| \leq~|a_{i_2i_1}(k_2+s-1)|(M^*+\varepsilon-\lambda^{s-1+k_2-k_1}\xi_1)
+(1-|a_{i_2i_1}(k_2+s-1)|)(M^*+\varepsilon)
\leq ~ M^*+\varepsilon-\lambda^{s+k_2-k_1}\xi_1, ~s=1,2,\dots.
$
We can further use the fact $k_2-k_1< T$ to obtain
\begin{align}\label{r2}
|x_{i_2}(k_1+s)| \leq M^*+\varepsilon-\lambda^{s}\xi_1, \ \ s=T, T+1,\dots.
\end{align}

We now reiterate the previous argument for the time interval $[k_1+T,k_1+2T)$. Again, there is a path from $i_1$ to any other node during the time interval $[k_1+T,k_1+2T)$. There exists a time instant $k_3\in [k_1+T+1,k_1+2T)$ such that either $i_1$ or $i_2$ is a neighbor of $i_3$ ($i_3$ is another node different from $i_1$ and $i_2$) at $k_3$. For any of the two cases we can deduce from (\ref{r1}) and (\ref{r2}) that for agent $i_3$, it must hold
$
|x_{i_3}(k_1+s)| \leq ~ M^*+\varepsilon-\lambda^{s}\xi_1,\ \ s=2T,2T+1,\dots.
$

The above analysis can be carried out to intervals $[k_1+2T,k_1+3T),\dots,[k_1+(n-2)T,k_1+(n-1)T)$, where $i_4,\dots,i_{n}$ can be found recursively until they include the whole network. We can therefore finally arrive at
$
\|x(k_1+(n-1)T)\|_{\infty} \leq ~ M^*+\varepsilon-\lambda^{(n-1)T}\xi_1
<~M^*-{\lambda^{(n-1)T}\xi_1}/{2},
$
for sufficient small $\varepsilon$ satifying $\varepsilon<{\lambda^{(n-1)T}\xi_1}/{2}$. Then, it follows from Lemma \ref{lem:invariant} that
\begin{align*}
\|x(k)\|_{\infty} <M^*-{\lambda^{(n-1)T}\xi_1}/{2},
\end{align*}
for all $k\geq k_1+(n-1)T$, which contradicts the fact that $\lim_{k\rightarrow \infty}\|x(k)\|_{\infty}=M^*$. Therefore,  $\lim_{k\rightarrow\infty}|x_i(k)|=M^*$, for all $i\in\mathcal{V}$.
\endproof

\subsection{Proof of Theorem \ref{thm1}: Existence of State Limits }\label{sec:directed2}

In this section, we show that $\lim_{k\rightarrow\infty}x_i(k)$ exists, for all $i\in\mathcal{V}$. Without loss of generality, $M^*$ is assumed to be nonzero and we fix any $i\in\mathcal{V}$.
Note that the fact that $\lim_{k\rightarrow\infty}|x_i(k)|=M^*$ include three possibilities: $\lim_{k\rightarrow\infty}x_i(k)=M^*$, $\lim_{k\rightarrow\infty}x_i(k)=-M^*$, or $x_i(k)$ switches between $-M^*$ and $M^*$ infinitely as $k\rightarrow\infty$. The last possibility actually means $\liminf_{k\rightarrow\infty}x_i(k)=-M^*$ and $\limsup_{k\rightarrow\infty}x_i(k)=M^*$.
We next prove the existence of the limit of $x_i(k)$ by showing that this last possibility cannot happen.

Suppose that we do have $\liminf_{k\rightarrow\infty}x_i(k)=-M^*$, and $\limsup_{k\rightarrow\infty}x_i(k)=M^*$. The following proof is based on a contradiction argument. We first use $M^*$ and $x_i(k)$ to bound the trajectory of $x_i$ after time instant $k$. Note that for all $s=1,2,\dots,$
\begin{align*}
x_{i}(k+s)=&~\sum_{j\in \mathcal{N}_{i}(k+s-1)}a_{ij}(k+s-1) x_j(k+s-1)
\\ \leq &~a_{ii}(k+s-1)x_{i}(k+s-1)
\\&+\!\!\!\!\!\!\!\!\sum_{j\in \mathcal{N}_{i}(k+s-1)\setminus\{ i\} }
\!\!\!\!\!\!|a_{ij}(k+s-1)|
|x_j(k+s-1)|
\\ \leq&~a_{ii}(k+s-1)x_{i}(k+s-1)
\\&+(1-a_{ii}(k+s-1))\|x(k)\|_{\infty}.
\end{align*}
It thus follows that
$
x_{i}(k+1)\leq \lambda x_{i}(k)+(1-\lambda)\|x(k)\|_{\infty}.$
Therefore, for all $k\geq \widehat{k}(\varepsilon)$, it follows that
$
x_{i}(k+1)\leq \lambda x_{i}(k)+(1-\lambda)(M^*+\varepsilon).
$
By a recursive analysis, we know that for all $k\geq \widehat{k}(\varepsilon)$ and all $s=1,2,\dots$,
\begin{align}
x_{i}(k+s)\leq&~\lambda^{s} x_{i}(k)+(1-\lambda^{s})(M^*+\varepsilon).\label{eq:x-bound}
\end{align}

Since $\liminf_{k\rightarrow\infty}x_i(k)=-M^*$, for any given $\varepsilon$, there exists an infinite sequence $\{\bar k_\chi\}_{\chi=0}^{\infty}$ such that $\bar k_{\chi}>\widehat{k}(\varepsilon)$ and $x_i(\bar k_\chi)\leq -(M^*-\varepsilon)$, $\chi = 0,1,\dots$. In addition, since $\limsup_{k\rightarrow\infty}x_i(k)=M^*$, for any $\bar k_\chi\in\{\bar k_\chi\}_{\chi=0}^{\infty}$, there exists a time instant $\bar k_{\overline{\chi}}>\bar k_{\chi}$ such that $x_i({\bar k}_{\overline{\chi}})\geq (M^*-\varepsilon)$.
By also noting that the state at each step is bounded by the previous step from \eqref{eq:x-bound},
there must exist a time instant $\bar k^*_{\chi}\in[\bar k_\chi,\bar k_{\overline{\chi}}]$ such that $-\lambda(M^*-\varepsilon)+(1-\lambda)(M^*+\varepsilon)\leq x_{i}(\bar k^*_{\chi})\leq~-\lambda^{2}(M^*-\varepsilon)+(1-\lambda^{2})(M^*+\varepsilon)$.

Therefore, for all $\bar k_\chi$, it follows that
\begin{align*}
|x_{i}(\bar k^*_{\chi})|\leq& ~\max\{|(1-2\lambda)M^*+\varepsilon|, |(1-2\lambda^2)M^*+\varepsilon|\}
\\ \leq &~ M^*-(1-\max\{|1-2\lambda|, |1-2\lambda^2|\})M^*+\varepsilon
\\ < &~ M^*-\frac{(1-\max\{|1-2\lambda|, |1-2\lambda^2|\})M^*}{2}
\end{align*}
if $\varepsilon$ is chosen sufficiently small as
\begin{displaymath}
\varepsilon<\frac{(1-\max\{|1-2\lambda|, |1-2\lambda^2|\})M^*}{2}.
\end{displaymath}
This contradicts the fact that $\lim_{k\rightarrow\infty}|x_i(k)|=M^*$ (which was shown in the beginning of Section \ref{sec:directed1}) and verifies that $\lim_{k\rightarrow\infty}x_i(k)$ exists. Therefore, we have proven Theorem \ref{thm1}.

\subsection{Proof of Theorem \ref{thm2}}\label{sec:undirected}

In this case, since $\mathcal{G}$ is infinitely jointly connected, the union graph $\mathcal{G}([k_1,\infty])$ is connected. We can therefore define
$
k_2:=\inf_k\big\{k\geq k_1, \mathcal{N}_{i_1}(k)\setminus\{i_1\}\neq \emptyset\big\}.
$
We denote $\mathcal{V}_1=\mathcal{N}_{i_1}(k_2)$. Obviously, we have that $|x_{i_1}(k_2)|=|x_{i_1}(k_1)|\leq M^*-\xi_1$, where $\xi_1=M^*-\alpha_{1}>0$ is defined as in \eqref{eq:i1}.
Therefore, following the similar analysis by which we obtained (\ref{r1}) and (\ref{r2}), we know that
$
|x_{i}(k_2+1)|\leq M^*+\varepsilon-\lambda\xi_1,\ i\in\mathcal{V}_1.
$

Similarly, since the union graph $\mathcal{G}([k_2+1,\infty])$ is connected, we can continue to define
$
k_3:=\inf_k\Big\{k\geq k_2+1: \bigcup_{i\in \mathcal{V}_1}\left(\mathcal{N}_{i}(k)\setminus\{i\}\right)\neq \emptyset\Big\}.
$
We also denote $\mathcal{V}_2=\bigcup_{i\in \mathcal{V}_1}\mathcal{N}_{i}(k_3)$. Note that $\{i_1\}\subseteq \mathcal{V}_1\subseteq \mathcal{V}_2$ with the definition of neighbor sets.
The fact that the graph is bidirectional guarantees that $k_3$ is not only the first time instant that there is an arc from ${\mathcal{V}_1}$ to another node, but also the first time instant that there is an arc from another node to ${\mathcal{V}_1}$ during the time interval $[k_2+1,k_3]$.
Therefore, we can apply Lemma \ref{lem:invariant} to the subset $\mathcal{V}_1$ for time interval $[k_2+1,k_3]$, and deduce that
$
 |x_{i}(k_3)|\leq M^*+\varepsilon-\lambda\xi_1,\  i\in{\mathcal{V}_1}.
$
It then follows from the same analysis that $
|x_{i}(k_3+1)|\leq M^*+\varepsilon-\lambda^2\xi_1,\ i\in\mathcal{V}_2.
$

The above argument can be carried out recursively for $\mathcal{V}_3$, $\mathcal{V}_4$, $\dots$ until $\mathcal{V}_m =\mathcal{V}$ for some constant $m\leq n-1$.  The corresponding $k_m$ can be found based on infinite joint connectedness condition, where $
|x_{i}(k_m+1)|\leq M^*+\varepsilon-\lambda^{m}\xi_1$, for all $ i\in \mathcal{V}$.
 This indicates that
\begin{align*}
\|x(k_m+1)\|_{\infty} \leq ~ M^*+\varepsilon-\lambda^{m}\xi_1<M^*-\lambda^{m}\xi_1/2,
\end{align*}
for sufficient small $\varepsilon$ satisfying $\varepsilon<\lambda^{n-1}\xi_1/2$. This contradicts the fact that
$\|x(k)\|_{\infty}\geq M^*-\varepsilon>M^*-\lambda^{m}\xi_1/2$, $\forall k\geq \widehat{k}(\varepsilon)$
(which was shown in the beginning of Section \ref{sec:directed1}). Therefore, it follows that $\lim_{k\rightarrow \infty}|x_i(k)|=M^*$, for all $i\in\mathcal{V}$.

Note that the proof given in Section \ref{sec:directed2} does not require connectivity. Therefore, using the same analysis as Section \ref{sec:directed2}, we can show that $\lim_{k\rightarrow \infty}x_i(k)$ exists, for all $i\in\mathcal{V}$. Therefore, we have proven Theorem \ref{thm2}.
\endproof

\subsection{Proof of Theorem \ref{thm3}}\label{sec:split}

It follows from Theorem \ref{thm1} that there exists a positive constant $M^*$ such that for all $i\in\mathcal{V}$, either $\lim_{k\rightarrow\infty}x_i(k)=M^*$ or $\lim_{k\rightarrow\infty}x_i(k)=-M^*$. We can therefore define two subsets of $\mathcal{V}$ as
$
\overline{\mathcal{V}}_1=\{i\in\mathcal{V}: \lim_{k\rightarrow\infty}x_i(k)=-M^*\},$
and
$
\overline{\mathcal{V}}_2=\{i\in\mathcal{V}: \lim_{k\rightarrow\infty}x_i(k)=M^*\}.$
Without loss of generality, we assume that $\overline{\mathcal{V}}_1$ is nonempty. Since the signed graph $\mathcal{G}^*=(V,\mathcal{E}^*)$ is sign consistent, the sign of each arc $(i,j)\in\mathcal{E}^*$ is denoted by $\varrho_{ij}$, where $\varrho_{ij}=+$ or $\varrho_{ij}=-$. We next show that if $(i,j)\in\mathcal{E}^*$ and $i,j\in\overline{\mathcal{V}}_1$, it is necessary that $\varrho_{ij}=+$.

Suppose it is not true, i.e., suppose $(i,j)\in\mathcal{E}^*$, $i,j\in\overline{\mathcal{V}}_1$, but $\varrho_{ij}=-$.
In the first place, it follows from the definitions of $\overline{\mathcal{V}}_1$ and $\overline{\mathcal{V}}_2$ that for any $\varepsilon>0$, there exists a positive constant $\widehat{k}_1(\varepsilon)$ such that for all $k\geq \widehat{k}_1(\varepsilon)$,
$
-M^*-\varepsilon \leq x_i\leq -M^*+\varepsilon,~ i\in \overline{\mathcal{V}}_1,
 M^*-\varepsilon \leq x_i\leq M^*+\varepsilon,~ i\in \overline{\mathcal{V}}_2.
$
Since the arc $(i,j)$ appears infinitely often in $\{\mathcal{G}_k\}_{0}^\infty$, it follows that there exists an infinite subsequence $\{\tilde{k}_\chi\}_{\chi=0}^{\infty}$ such that $\tilde{k}_\chi\geq \widehat{k}_1(\varepsilon)$ and $(i,j)\in \mathcal{E}_{\tilde{k}_\chi}$ for all $\chi=0,1,\dots$. Consider any $\tilde{k}_\chi$. We know that
\begin{align*}
x_{j}(\tilde{k}_\chi+1)= &~a_{jj}(\tilde{k}_\chi)x_{j}(\tilde{k}_\chi)+a_{ji}(\tilde{k}_\chi)x_{i}(\tilde{k}_\chi)
\\&+\!\!\!\!\!\!\!\!\sum_{l\in \mathcal{N}_{j}(\tilde{k}_\chi)\setminus\{ j,i\} }
\!\!\!\!\!\!a_{jl}(\tilde{k}_\chi)
x_l(\tilde{k}_\chi)
\\ \geq&~~-a_{jj}(\tilde{k}_\chi)(M^*+\varepsilon)+a_{ji}(\tilde{k}_\chi)(-M^*+\varepsilon)
\\& ~-\sum_{l\in \mathcal{N}_{j}(\tilde{k}_\chi)\setminus\{ j,i\} }|a_{jl}(\tilde{k}_\chi)|(M^*+\varepsilon)
\\=&~-a_{jj}(\tilde{k}_\chi)(M^*+\varepsilon)+a_{ji}(\tilde{k}_\chi)(-M^*+\varepsilon)
\\&-(1-a_{jj}(\tilde{k}_\chi)-|a_{ji}(\tilde{k}_\chi)| )(M^*+\varepsilon)
\\ \geq&~ (2\lambda-1)M^*-\varepsilon,
\end{align*}
where we have used the fact that $a_{ji}(\tilde{k}_\chi)\leq 0$ from Assumption \ref{assm1}. If we choose $\varepsilon$ sufficiently small satisfying $\varepsilon<\lambda M^*$, it then follows that $x_{j}(\tilde{k}_\chi+1) >  -M^*+\varepsilon$. Note that  $x_{j}(\tilde{k}_\chi+1)>-M^*+\varepsilon$ holds for all $\chi=0,1,\dots$. This shows that $\liminf_{\chi \rightarrow\infty}x_{j}(\tilde{k}_\chi+1)> -M^*+\varepsilon$, which contradicts that $j\in \overline{\mathcal{V}}_1$.
We thus know that if $(i,j)\in\mathcal{E}^*$ and $i,j\in\overline{\mathcal{V}}_1$, then $\varrho_{ij}=+$.

Next, since $\mathcal{G}^*=(\mathcal{V},\mathcal{E}^*)$ is structurally balanced,
$\mathcal{V}$ is divided into two disjoint nonempty subsets $\mathcal{V}_1$ and $\mathcal{V}_2$, where $\varrho_{kl}=+$, for all $k,l\in\mathcal{V}_1$ and $k,l\in\mathcal{V}_2$, and $\varrho_{kl}=-$ for all $k\in\mathcal{V}_1$, $l\in\mathcal{V}_2$ and $k\in\mathcal{V}_2$, $l\in\mathcal{V}_1$.

It therefore follows that $\overline{\mathcal{V}}_1 \subseteq \mathcal{V}_1$ from the reasoning that $(i,j)\in\mathcal{E}^*$ and $i,j\in\overline{\mathcal{V}}_1$ implies $\varrho_{ij}=+$. Therefore, we know that $\overline{\mathcal{V}}_2$ is nonempty. In addition, $\varrho_{kl}=+$, for all $k,l\in\mathcal{V}_2$. We thus know that $\overline{\mathcal{V}}_1 =\mathcal{V}_1$ and $\overline{\mathcal{V}}_2 =\mathcal{V}_2$.
%Moreover, for the special of $M^*=0$, we know that $\overline{\mathcal{V}}_2$ is empty. But we have shown that %$\overline{\mathcal{V}}_2$ is nonempty. Thus, this case cannot hold.
Therefore, it follows that $\lim_{k\rightarrow\infty}x_i(k)=M^*$, $i\in\mathcal{V}_1$ and $\lim_{k\rightarrow\infty}x_i(k)=-M^*$, $i\in\mathcal{V}_2$ for every initial state $x(0)\in \mathbb{R}^{n}$ , where $M^*$ is a constant.

%\subsection{Proof of Theorem \ref{thm4}}
%
%Note that the proof given in Section \ref{sec:split} does not rely on connectivity. Therefore, the proof of Theorem \ref{thm4} is exactly the same as that of Theorem \ref{thm3}.

\section{Numerical Example}\label{sec:simulation}

Consider the following discrete-time Kuramoto oscillator system with antagonistic and cooperative links:
\begin{equation}
\theta_i(k+1)=\theta_i(k)-\mu\sum_{j\in \mathcal{N}_i(k)\setminus \{i\}}\sin\Big(\theta_i(k)-R_{ij}(k)\theta_j(k)\Big),\label{eq:oscillator}
\end{equation}
where $\theta_i(k)$ denotes the state of node $i$ at time $k$, $\mu>0$ is the stepsize, and $R_{ij}(k)\in\{1,-1\}$ represents the cooperative or antagonistic
relationship between node $i$ and node $j$. Note that with $R_{ij}(k)\equiv 1$, system \eqref{eq:oscillator} corresponds to the classical Kuramoto oscillator model (\cite{kuramoto1}). Let $\delta\in (0,\frac{\pi}{2})$ be a given constant and suppose $\theta_i(0)\in(-\frac{\pi}{2}+\delta,\frac{\pi}{2}-\delta)$ for all $i\in\mathcal{V}$. Here $\delta$ can be any positive constant sufficiently small.
System \eqref{eq:oscillator} can be rewritten as
\begin{align*}
\theta_i(k+1)=&~\theta_i(k)-\mu\!\!\!\sum_{j\in \mathcal{N}_i(k)\setminus \{i\}}\!\!\!\frac{\sin(\theta_i(k)-R_{ij}(k)\theta_j(k))}{\theta_i(k)-R_{ij}(k)\theta_j(k)}
\\&~\times(\theta_i(k)-R_{ij}(k)\theta_j(k)).
\end{align*}
Note that the function
$\sin x / x$
is well-defined for $x\in(-\infty,\infty)$.
Therefore, we can define
\begin{displaymath}
a_{ij}(\theta,k)=\\ \frac{\sin(\theta_i(k)-R_{ij}(k)\theta_j(k))}{\theta_i(k)-R_{ij}(k)\theta_j(k)}R_{ij}(k), \ \ j\in \mathcal{N}_i(k)\setminus \{i\},
\end{displaymath}
and $a_{ii}(\theta,k)=1-\mu\sum_{j\in \mathcal{N}_i(k)\setminus \{i\}}|a_{ij}(\theta,k)|$, where $\theta=(\theta_1,\theta_2,\dots,\theta_n)\T$, so that \eqref{eq:oscillator} is re-written into the form of \eqref{eq:origin}.

Lemma \ref{lem:invariant} ensures that
$
0<\lambda^*\leq  \frac{\sin(\theta_i(k)-R_{ij}(k)\theta_j(k))}{\theta_i(k)-R_{ij}(k)\theta_j(k)}\leq 1,
$
where $\lambda^*=\frac{\sin(\pi-2\delta)}{\pi-2\delta}$. This gives us $|a_{ij}(\theta,k)|\geq \lambda^*$, for all $i$ and $j\in \mathcal{N}_i(k)\setminus \{i\}$. In addition,
by selecting $\mu<\frac{1-\lambda^*}{n}$, we can guarantee that $|a_{ii}(\theta,k)|\geq \lambda^*$ for all $\theta$ and $k$.
Therefore, we can use Theorems \ref{thm1}, \ref{thm2}, and \ref{thm3} to study the behavior of Kuramoto oscillator with antagonistic links.

%Therefore, given $\mu<\frac{1-\lambda^*}{n}$, it follows that $\lim_{k\rightarrow\infty}(|\theta_i(k)|-|\theta_j(k)|)=0$, $i,j\in \mathcal{V}$ under \eqref{eq:oscillator}, if $\{\mathcal{G}_k\}_0^\infty$ is uniformly jointly strongly connected for unidirectional graphs, or infinitely jointly connected for bidirectional graphs according to Theorems \ref{thm1} and \ref{thm2}.

\begin{figure}[t]
\begin{minipage}[t]{0.3\linewidth}
~~~~~~~~~~ \xymatrix{
1\ar@{->}[r]^-{-} &  3\\
&  2
}
\end{minipage}
\begin{minipage}[t]{0.3\linewidth}
\centering
~~ \xymatrix{
1 &  3\ar@{->}[d]^-{-}\\
&  2
}
\end{minipage}
\begin{minipage}[t]{0.3\linewidth}
\centering
~~  \xymatrix{
1\ar@{<-}[rd]_-{+}  &  3\ar@{<-}[d]_-{+}\\
&  2
}
\end{minipage}\caption{Graphs $\mathcal{G}_1$, $\mathcal{G}_2$,  $\mathcal{G}_3$ considered in the example.}\label{fig-g123}
\end{figure}
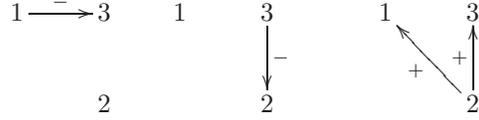

We next verify the theoretical results using simulations. For the case of unidirectional topology, we assume that
the topology switches periodically as $\begin{tabular} {cc}
 \xymatrix{
\mathcal{G}_1 \ar@{->}[r] &  \mathcal{G}_2\ar@{->}[r] &  \mathcal{G}_3\ar@{->}[r]&  \mathcal{G}_1\ar@{->}[r] &\dots
} \end{tabular}$
at time instants $\eta_l=l~s$, $l=1,2,\dots$, where
$\mathcal{G}_1$, $\mathcal{G}_2$, $\mathcal{G}_3$ are represented in Fig.~\ref{fig-g123}.
The system matrices associated with $\mathcal{G}_1$, $\mathcal{G}_2$, $\mathcal{G}_3$ are given by
$
A_1=\left[\begin{smallmatrix} 1&0&0\\0 &1& 0\\-0.5 &0 &0.5
\end{smallmatrix}\right], ~A_2=\left[\begin{smallmatrix} 1&0&0\\0 &0.5& -0.5\\0 &0 &1
\end{smallmatrix}\right],A_3=\left[\begin{smallmatrix} 0.5&0.5&0\\0 &1& 0\\0 &0.5 &0.5
\end{smallmatrix}\right].$

\begin{figure}[t]
\begin{center}
\includegraphics[width=6.5cm,height=4.5cm]{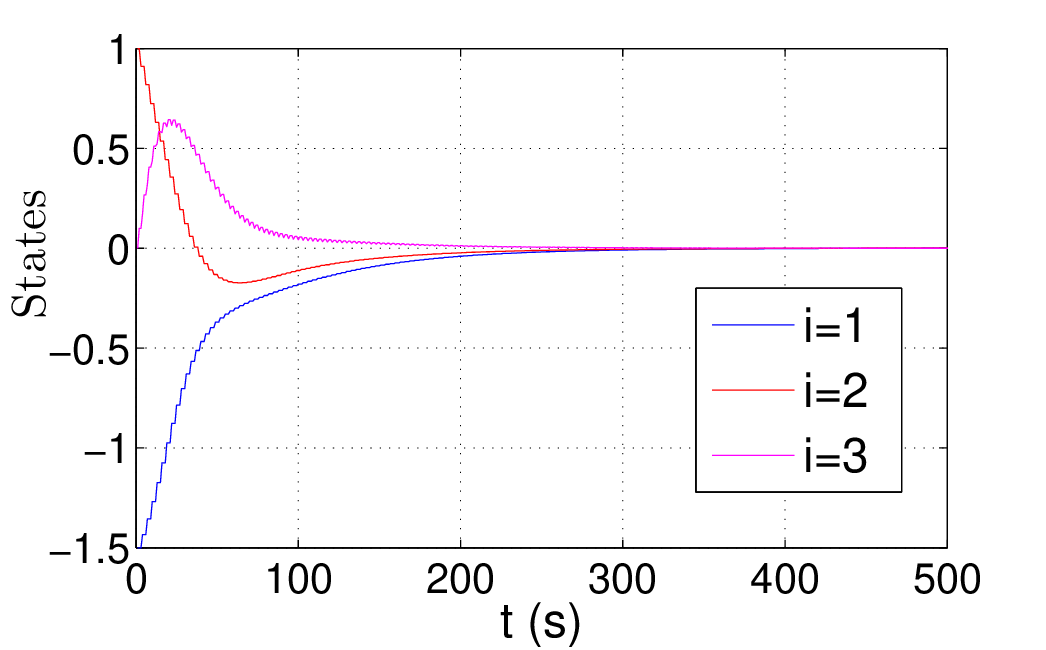}
\caption{Convergence for unidirectional topology}\label{fig-directed}
 \end{center}
\end{figure}

The initial state is $x(0)=(-1.5,1,0)\T$ and $\mu=0.1$. Fig.~\ref{fig-directed} shows the convergence of states over unidirectional switching topologies. We see that the absolute values of the states converge for this group of oscillators
with antagonistic interactions and switching topologies, in accordance with the conclusion from Theorem~\ref{thm1}. Note all agent states converge to zero, instead of achieving bipartite consensus.

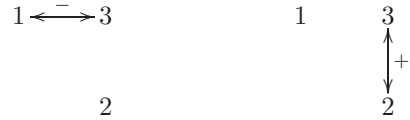
\begin{figure}[t]
\begin{minipage}[t]{0.45\linewidth}
\centering
~~ \xymatrix{
1 &  3\ar@{<->}[l]_-{-}\\
&  2
}
\end{minipage}
\begin{minipage}[t]{0.45\linewidth}
\centering
~ \xymatrix{
1 &  3\\
&  2\ar@{<->}[u]_-{+}
}
\end{minipage}\caption{Graphs $\mathcal{G}_4$ and $\mathcal{G}_5$ considered in the example.}\label{fig-g45}
\end{figure}

\begin{figure}[t]
\begin{center}
\includegraphics[width=6.5cm,height=4.5cm]{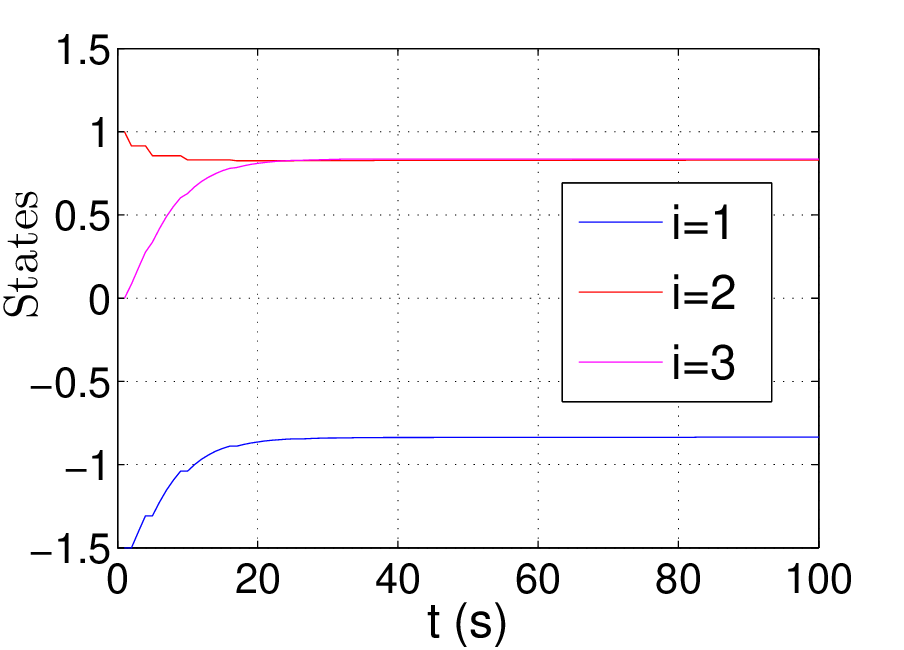}
\caption{Convergence for bidirectional topology}\label{fig-undirected}
 \end{center}
\end{figure}

For the case of bidirectional topology, we assume the topology switches between $\mathcal{G}_4$ and $\mathcal{G}_5$ in Fig. \ref{fig-g45}. The topology is $\mathcal{G}_4$ except at time intervals $[l^2,l^2+1]$, where the topology is $\mathcal{G}_5$, $l=1,2,\dots$. The signed matrices associated with $\mathcal{G}_4$, $\mathcal{G}_5$ are
$
A_4=\left[\begin{smallmatrix} 0.5&0&-0.5\\0 &1& 0\\-0.5 &0 &0.5
\end{smallmatrix}\right], ~A_5=\left[\begin{smallmatrix} 1&0&0\\0 &0.5& 0.5\\0 &0.5 &0.5
\end{smallmatrix}\right].
$
The initial states and $\mu$ are the same as previously. Fig.~\ref{fig-undirected} shows the convergence of states over this bidirectional switching topology. We see that
the absolute values of the agent states converge to the same limit, in accordance with the conclusion from Theorems \ref{thm2} and \ref{thm3}. Bipartite consensus is achieved in this case.

\section{Conclusions}\label{sec:conclusion}

In this paper, we
studied the consensus problem of multi-agent systems over
cooperative-antagonistic networks in a discrete-time setting. Both unidirectional and bidirectional topologies were considered. It was proven that the limits of all agent states exist and reach a consensus in absolute values if the topology is uniformly jointly strongly connected or infinitely jointly connected. We also gave an example to show that uniform quasi-strong connectedness is not sufficient to guarantee consensus in absolute values. We further proposed sufficient conditions for bipartite consensus to
be achieved over networks with joint connectivity. Examples were given to explain coordination of multiple nonlinear systems with antagonistic interactions using the proposed algorithms. Future works include investigating time-delay influence and other types of antagonistic interaction models.

%considering antagonistic interactions for consensus problems for higher-order dynamics dynamics and

%%%%%%%%%%%%%% begin Bibliography %%%%%%%%%%%%%%%%%

%\bibliographystyle{IEEEtran}
\bibliographystyle{elsarticle-harv}
\bibliography{refs}
%%%%%%%%%%%%%%%% end Bibliography %%%%%%%%%%%%%%%%%

%%%%%%%%%%%%%%%% biosketches and photos %%%%%%%%%%%%%%%%%

\end{document}